\documentclass{amsart}

\usepackage{graphicx}
\usepackage{amsmath}
\usepackage{amsfonts,amssymb}
\usepackage[latin1]{inputenc}
\usepackage{amsthm}
\usepackage{verbatim}
\usepackage{array}
\usepackage{graphicx,color}
\usepackage{algorithm}
\usepackage[noend]{algpseudocode}

\makeatletter
\def\BState{\State\hskip-\ALG@thistlm}
\makeatother

\newtheorem{problem}{Problem}

\begin{document}

\title{Numerical calculation of extremal Steklov eigenvalues in 3D and 4D}
\author{Pedro R. S. Antunes}

\address{Univ. Aberta, Dep. of Sciences and Technology,
  Rua da Escola Polit\'ecnica 141-7, P-1269-001 Lisboa, Portugal, and
  Univ. Lisboa, Faculty of Sciences, Group of Mathematical Physics,
Edif\'icio 6, Piso 1,  Campo Grande,
  P-1749-016 Lisboa, Portugal.} \email{prantunes@fc.ul.pt}


\date{\today}

\begin{abstract}
We develop a numerical method for solving shape optimization of functionals involving Steklov eigenvalues and apply it to the problem of maximization of the $k$-th Steklov eigenvalue, under volume constraint. A similar study in the planar case was addressed in [E. Akhmetgaliyev, C.-Y. Kao and B. Osting, SIAM J. Control Optim. \textbf{55}(2), 1226-1240, (2017)] using the boundary integral equation method. Here we extend that study to the 3D and 4D cases, using the Method of Fundamental Solutions as forward solver. 
\end{abstract}

\maketitle

\pagestyle{myheadings} \thispagestyle{plain} \markboth{Pedro R. S.
Antunes}{Optimization of Steklov eigenvalues}

\section{Introduction\label{sec:intro}}
Let $\Omega\subset\mathbb{R}^d$ be a bounded open set
and consider the second order Steklov eigenvalue problem,
\begin{equation}\label{eigprobdir}
\left\{\begin{array}{ll} \Delta w =0 & \textrm{in } \Omega\\[1mm]
\frac{\partial w}{\partial n}=\sigma w & \textrm{on } \partial\Omega,\end{array}\right.
\end{equation}
\noindent  
defined in $H^1(\Omega)$. We will
denote the eigenvalues by
$0=\sigma_0(\Omega)\leq\sigma_1(\Omega)\leq\sigma_2(\Omega)\leq...\rightarrow\infty$ where each eigenvalue $\sigma_k(\Omega)$
is counted with
its multiplicity and the corresponding orthonormal real
eigenfunctions by $w_i,\ i=0,1,2,...$. The eigenvalues can also be defined through a variational characterization
\begin{equation}
\label{carvarstek}
\sigma_k(\Omega)=\min_{S\in\mathcal{S}_{k+1}}\max_{w\in S\backslash\left\{0\right\}}\frac{\int_\Omega|\nabla w|^2dx}{\int_{\partial\Omega}w^2ds_x},
\end{equation}
where $\mathcal{S}_{k+1}$ is the family of all subspaces of dimension $k+1$ in $H^1(\Omega)$ (eg.~\cite{BBG}).

Some shape optimization problems for Steklov eigenvalues were already considered in the literature. For example, in~\cite{weinst} the author proved that the disk maximizes the first non trivial eigenvalue of the problem
\begin{equation}\label{eigprobdirweins}
\left\{\begin{array}{ll} \Delta w =0 & \textrm{in } \Omega\\[1mm]
\frac{\partial w}{\partial n}=\sigma\rho w & \textrm{on } \partial\Omega,\end{array}\right.
\end{equation}
\noindent  
among simply connected planar domains with a fixed mass $M(\Omega)=\int_{\partial\Omega}\rho ds_x$, where $\rho\geq0$ is a $L^\infty(\partial\Omega)$ function on the boundary. The extension of this result to non simply connected domains is an open problem. For higher eigenvalues, in~\cite{HPS} the authors proved a result that implies the following inequality for simply connected planar domains,
\begin{equation}
\label{hersh}
\sup\left\{\sigma_k(\Omega)M(\Omega):\ \Omega\subset\mathbb{R}^2\right\}\leq 2\pi k, k\in\mathbb{N}.
\end{equation}
\noindent
In~\cite{GP}, it was proved that the previous inequality is sharp and attained by a sequence of simply connected domains degenerating into a disjoint union of $k$ disks of the same area.

The connection of extremal Steklov eigenvalue problems and the problem of generating free boundary minimal surfaces in the Euclidean ball was studied in \cite{Fraser} and this connection was explored numerically in \cite{oudet_ko}.\\

In this work, we will consider the following shape optimization problem
\begin{equation}
\label{shopteigvalprob}
\sigma_k^{\star}=\max_{\Omega\subset\mathbb{R}^d}\left\{\sigma_k(\Omega),\ |\Omega|=1\right\},\ k=1,2,....
\end{equation}
where $\sigma_k$ are the eigenvalues of problem \eqref{eigprobdir}
and taking into account the properties
\[\sigma_k(t\Omega)=\frac{1}{t}\sigma_k(\Omega),\ \forall t>0\]
and
\[|t\Omega|=t^d|\Omega|,\ \forall t>0\]
the problem \eqref{shopteigvalprob} is equivalent to

\begin{problem}
\label{shopteigvalprobnorm}
Given $k=1,2,...$, determine
\[\sigma_k^{\star}=\max_{\Omega\subset\mathbb{R}^d}\left\{\sigma_k(\Omega)|\Omega|^{1/d}\right\}.\]
\end{problem}
The existence of an optimal set was proved in~\cite{BBG} and some numerical studies for the optimal solutions of this problem in 2D can be found in~\cite{ak-kao-osting,Bogosel,BBG}. As pointed out in~\cite{ak-kao-osting}, the optimal domains for this problem in 2D are well structured. In particular, the optimizer for $\sigma_k$ seems to have $k$-fold symmetry and the corresponding optimal eigenvalue has multiplicity 2, if $k$ is even and multiplicity 3, if $k\geq3$ is odd. Also some progress on the analytical study of this problem has been made recently. In~\cite{osting2018}, using a perturbation argument the authors proved that the disk is not the maximizer for higher even numbered Steklov eigenvalues. In this work we will address a numerical study for the optimizers in 3D and 4D.\\

\section{Shape derivatives}
\label{shapeder}
We will consider the numerical solution of Problem~\ref{shopteigvalprobnorm} using a gradient-type method. In this context it is convenient to make use of the formulas for the shape derivatives of Steklov. Consider an application $\Psi:t\in[0,T[\rightarrow W^{1,\infty}(\mathbb{R}^d,\mathbb{R}^d)$ for which $\Psi(t)=I+tV$, where $W^{1,\infty}(\mathbb{R}^d,\mathbb{R}^d)$ is the set of
bounded Lipschitz maps from $\mathbb{R}^d$ into itself, $I$ is the
identity and $V$ is a given deformation field. 

We will use the notation 
$\Omega_t=\Psi(t)(\Omega)$ and $\sigma(t):=\sigma(\Omega_t)$ and we assume that $\sigma(0)$ is simple.

Define the function $\mathcal{V}(t)=|\Omega_t|$, which gives the volume of the domain $\Omega_t$. Then the shape derivative of $\mathcal{V}$ is given by
\begin{equation}
\label{dervolume} \mathcal{V}'(0)=\int_{\partial\Omega} V.n\, d s_x.
\end{equation}

The shape derivative for a Steklov eigenvalue, which can be found in~\cite{DKL,ak-kao-osting}, is given by
\[\sigma'(0)=\int_{\partial\Omega}\left(|\nabla w|^2-2\sigma^2w^2-\sigma\mathcal{H} w^2\right)V.nds_x,\]
where $\mathcal{H}$ is the mean curvature.

\section{Numerical methods}
\label{nummethods}

In this section we describe briefly the numerical methods that were used to solve Problem~\ref{shopteigvalprobnorm}.

\subsection{Parametrization of the domains}
For the numerical solution of the Problem~\ref{shopteigvalprobnorm}, we assume that $\Omega$ is a star-shaped domain whose boundary can be parameterized by
\begin{equation}\label{3Dsphere}
\left\{\begin{array}{l} x=r(\theta,\phi)\sin(\theta)\cos(\phi) \\[1mm]
y=r(\theta,\phi)\sin(\theta)\sin(\phi)\\[1mm]
z=r(\theta,\phi)\cos(\theta)
 \end{array}\right.
\end{equation}
with $r(\theta,\phi)>0$ for $\theta\in[0,\pi]$ and $\phi\in[0,2\pi[$, and

\begin{equation}\label{4Ddomain}
\left\{\begin{array}{l} x=r(\beta,\theta,\phi)\sin(\beta)\sin(\theta)\cos(\phi) \\[1mm]
y=r(\beta,\theta,\phi)\sin(\beta)\sin(\theta)\sin(\phi)\\[1mm]
z=r(\beta,\theta,\phi)\sin(\beta)\cos(\theta)\\[1mm]
w=r(\beta,\theta,\phi)\cos(\beta)\\[1mm]
 \end{array}\right.
\end{equation}
where $r(\beta,\theta,\phi)>0$, for $\beta\in[0,\pi]$, $\theta\in[0,\pi]$ and $\phi\in[0,2\pi[$, respectively for 3D and 4D domains.

We define the family of 3D spherical harmonics

\[S_l^m(\theta,\phi)=\left\{ \begin{array}{ll}
\sqrt{2}k_l^m\cos(m\phi)P_l^m(\cos(\theta)) & \mbox{ if } m>0, \\[5pt]
k_l^0P_l^0(\cos(\theta)) & \mbox{ if } m=0, \\[5pt]
\sqrt{2}k_l^m\sin(-m\phi)P_l^{-m}(\cos(\theta)) & \mbox{ if } m<0,\end{array}\right.\]
where $P_l^m$ is an associated Legendre polynomial and
\[k_l^m=\sqrt{\frac{(2l+1)(l-|m|)!}{4\pi(l+|m|)!}}.\]
The 4D hyper-spherical harmonics are defined by
\begin{equation}
\label{spharm}
S_{nl}^{m}(\beta,\theta,\phi)=c_{n,l,m}\sin^l(\beta)C_{n-l}^{l+1}\left(\cos(\beta)\right)
S_l^m(\theta,\phi)(\beta,\theta,\phi),\ \begin{array}{l} n=0,1,2,...\\[1mm]
0\leq l\leq n\\[1mm]
-l\leq m\leq l\end{array}, \begin{array}{l} \beta\in[0,\pi]\\[1mm]
\theta\in[0,\pi]\\[1mm]
\phi\in[0,2\pi[\end{array},
\end{equation}
where $S_l^m$ are 3D spherical harmonics, $C_{n-l}^{l+1}$ are Gegenbauer polynomials and
\[c_{n,l,m}=2^{l+\frac{1}{2}}\sqrt{\frac{(n+1)\Gamma(n-l+1)}{\pi\Gamma(n+l+2)}}\Gamma(l+1).\]

\noindent
The function $r$ is expanded in 3D and 4D respectively in terms of 3D spherical harmonics, for a fixed $N\in\mathbb{N}$,
\begin{equation}
\label{linqcomb3d}
r(\theta,\phi)=\sum_{l=0}^N\sum_{m=-l}^la_{l,m}S_l^m(\theta,\phi),
\end{equation}
and 4D hyper-spherical harmonics
\begin{equation}
\label{rad}
r(\beta,\theta,\phi)=\sum_{n=0}^{N}\sum_{l=0}^n\sum_{m=-l}^{l}a_{n,l,m}S_{nl}^{m}(\beta,\theta,\phi)
\end{equation}
and the optimization is performed by searching for optimal coefficients in these expansions.

\subsection{Generation of points for the Method of Fundamental Solutions}
In~\cite{AO} it was proposed a fast algorithm for the generation of points on the boundary of 3D and 4D star-shaped domain. It is based on an almost uniform distribution of points on the unitary sphere. Another numerical approach to obtain a quasi-equidistant point distribution over the surface of a sphere was proposed in~\cite{Serranho}. The problem of distributing some points on a sphere in order to maximize the minimum distance between pairs of points is known as Tammes problem (\cite{tammes,conway}) and the optimal distribution was already found for some particular small numbers $N$ (eg.~\cite{fejes,danzer,musin}). 
 
\begin{figure}[ht]
\centering \includegraphics[width=0.48\textwidth]{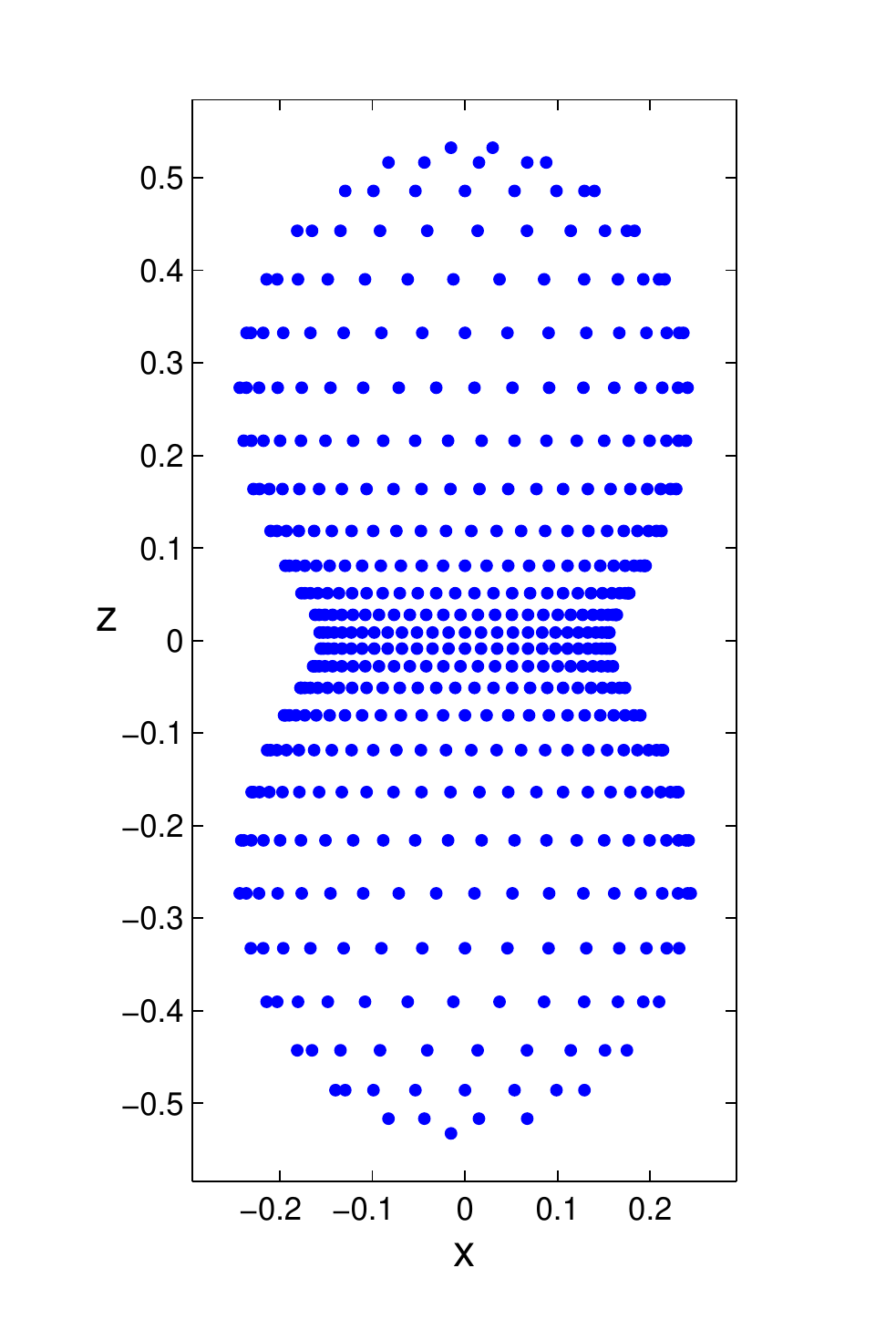}
\centering \includegraphics[width=0.48\textwidth]{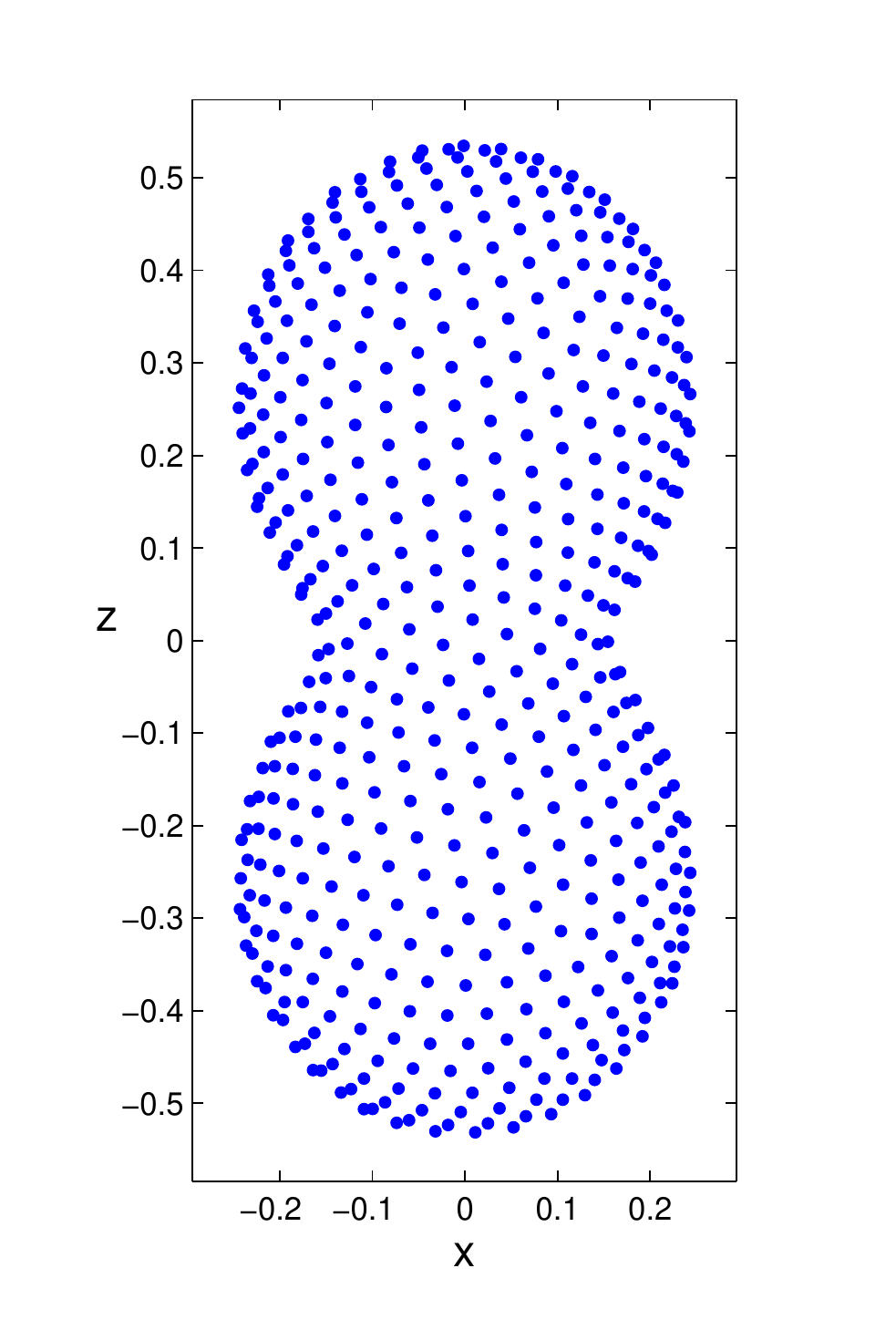}
\caption{Two distribution of points on the boundary of a star-shaped domain.} \label{fig:points}
\end{figure}

The location of points on the boundary of a general star-shaped domain can be obtained directly from \eqref{4Ddomain}, mapping a sample of points almost uniformly distributed on the surface of the sphere to points on the boundary of the domain. However, we have the effect of the function $r$ and in some cases this can generate clusters of nodes (see Figure~\ref{fig:points}-left).
 
Another strategy is to calculate the points $x_i$, $i=1,...,N$ that minimize the Riesz energy,
\[E(N,s)=\sum_{1\leq i<j\leq N}\frac{1}{|x_i-x_j|^s}.\]
The Coulomb potential that models electrons repelling each other is the case $s=1$. To illustrate the results obtained with both numerical algorithms we considered a 3D star-shaped domain parametrized by \eqref{3Dsphere}, where
\[r(\theta,\phi)=1+0.4S_2^0(\theta,\phi).\]
In Figure~\ref{fig:points}-left we plot 1000 points obtained with the distribution of points considered in~\cite{AO}. Then, this distribution of points was improved by minimizing the Riesz energy, with $s=3$ and we obtained the distribution of nodes of Figure~\ref{fig:points}-right. This approach is much more expensive from the computational point of view, but provides a much better distribution of nodes. We considered an analogous algorithm for placing the nodes on the boundary of 4D domains.

We will follow the choice for the source points of the MFS proposed in~\cite{AlvAnt,AA}. We take $M_\mathcal{C}$ collocation points $x_i,\ i=1,...,M_\mathcal{C}$ almost uniformly distributed on the boundary of the domain obtained by using the algorithm described above and for each of these points we calculate the outward unitary vector $n_i$, which is normal to the boundary at $x_i$. The source points are defined by
\[y_i=x_i+\delta\ n_i,\]
where $\delta$ is a parameter chosen such that the source points remain outside $\bar{\Omega}$.

\subsection{Numerical solution of the eigenvalue problems}
The Steklov eigenvalue Problem~\ref{eigprobdir} will be solved by the Method of Fundamental Solutions (MFS). We take the fundamental solution of the Laplace equation in $\mathbb{R}^d$, $d\geq3$ (e.g.~\cite{evans}),
\begin{equation}
\label{solfund}
\Phi(x)=\frac{1}{d(d-2)\alpha(d)|x|^{d-2}},
\end{equation}
where $\alpha(d)$ denotes the volume of the unit ball in $\mathbb{R}^d$.

The MFS approximation is a linear combination
\begin{equation}
\label{mfsapp}
u(x)\approx\tilde{u}(x)=\sum_{j=1}^M\beta_j\phi_j(x),
\end{equation}
where
\begin{equation}
\label{ps} \phi_{j}=\Phi(\cdot-y_{j})
\end{equation}
are $M$ point sources centered at some points $y_{j}$ that are
placed on an admissible source set $\hat{\Gamma}$, which is assumed to be the boundary of a bounded open set $\hat{\Omega}$ such that $\bar{\Omega}\subset\hat{\Omega}$, with $\hat{\Gamma}$ surrounding $\partial\Omega$. The MFS approximation \eqref{mfsapp} can be seen as a discretization of the single layer operator
\begin{align*}
S : &C(\hat{\Gamma})\rightarrow C(\partial\Omega)\\
(S\varphi)(x)&=\int_{\hat{\Gamma}}\Phi(x-y)\varphi(y)ds_y
\end{align*}
which is known to be injective, for $d>2$ (e.g.~\cite{Kress}).
By construction, the MFS approximation satisfies the Laplace equation because it is built using shifts of the fundamental solution. Moreover, since it is a mesh and integration free method it is particularly suitable to deal with boundary value problems defined in 3D or 4D domains, as those considered in this paper. For more details about the MFS, we refer to the following works \cite{Kupradze,Bog,FK,AlvAnt,bb,AA,Ant,Bogosel2,Ant2}.

The approximation of the Steklov eigenvalues can be calculated by solving generalized matrix eigenvalue problems. We define the matrices $\mathbf{A}$ and $\mathbf{B}$, where
\[ (\mathbf{A})_{i,j}=\partial_{n_i}\Phi(x_i-y_j)\quad(\mathbf{B})_{i,j}=\Phi(x_i-y_j).\]

The eigenvalues are calculated by solving the generalized matrix eigenvalue problem
\[\mathbf{A}X=\Lambda\mathbf{B}X\]
using the Matlab routine \texttt{eigs}.

\begin{figure}[ht]
\includegraphics[width=0.48\textwidth]{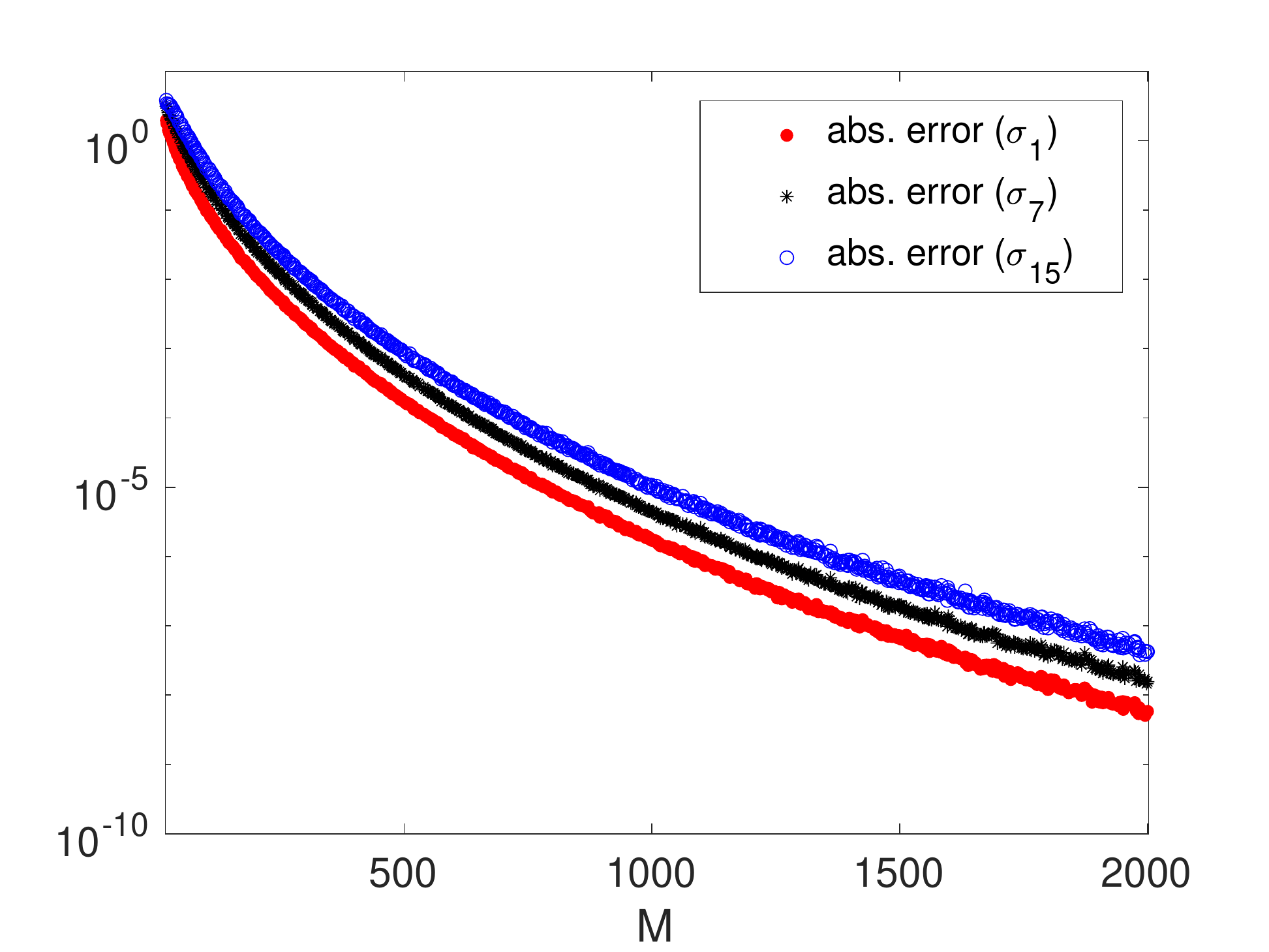}
\includegraphics[width=0.48\textwidth]{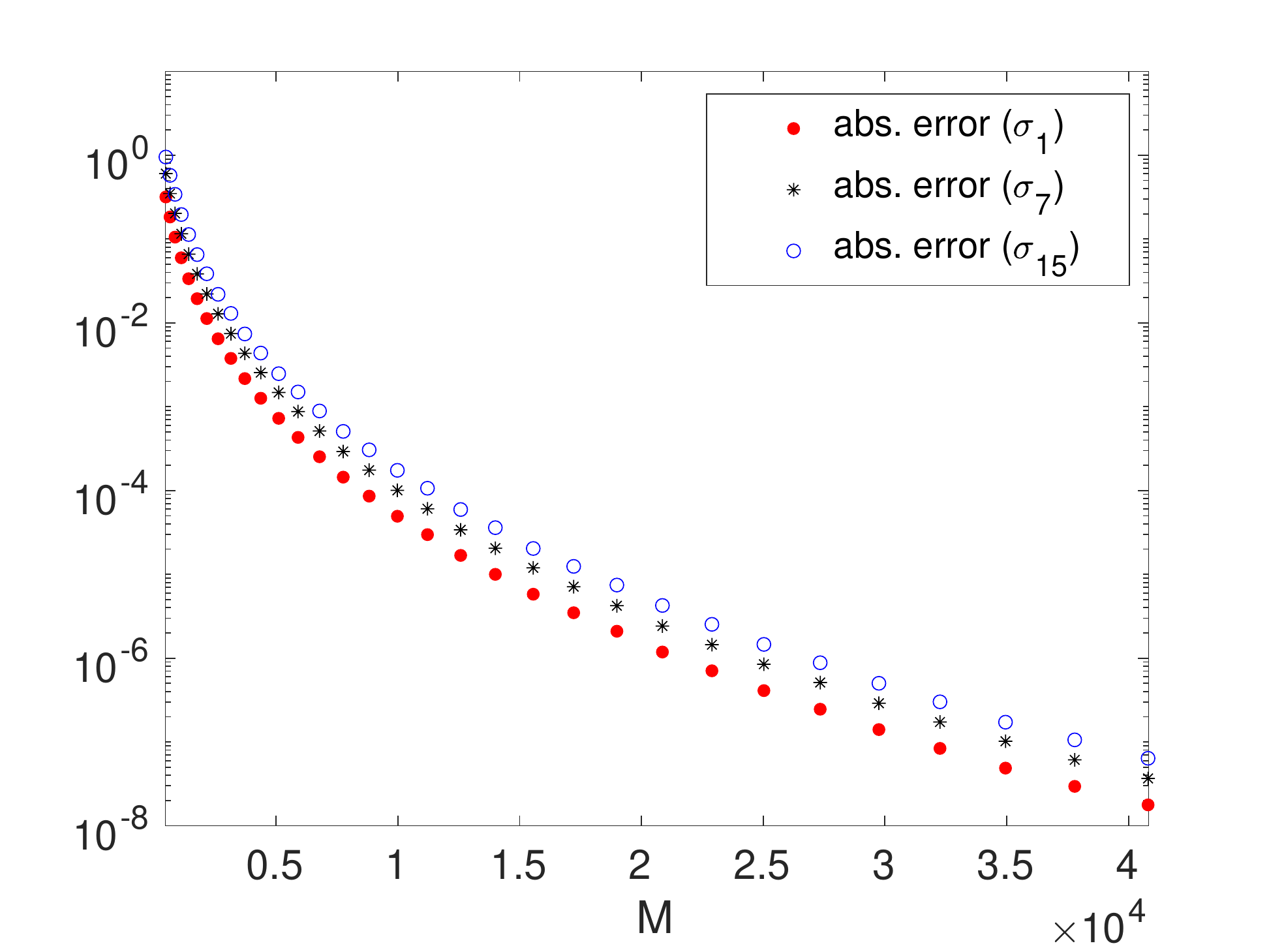}
\caption{Convergence curve for the MFS approximation of $\sigma_k$, $k=1,7,15$ of the ball of unit volume in 3D (left plot) and 4D (right plot).} \label{fig:conv}
\end{figure}

Next, we test our algorithm for the calculation of Steklov eigenvalues in the case of 3D and 4D balls, for which we know the exact solutions. Figure~\ref{fig:conv} shows the absolute error of the approximations for three eigenvalues, $\sigma_1$, $\sigma_7$ and $\sigma_{15}$ in 3D (left plot) and 4D (right plot), which were obtained for $\delta=0.2.$

\subsection{Optimization algorithm}

We define $\mathcal{V}\in\mathbb{R}^P$ to be the vector of all the coefficients in the linear combinations \eqref{linqcomb3d} and \eqref{rad}. Then, we define the cost function
\[\mathcal{C}_k(\mathcal{V}):=\sigma_k(\Omega)|\Omega|^{\frac{1}{d}},\]
where $\Omega$ is the domain whose boundary is obtained from $\mathcal{V}$, by using \eqref{linqcomb3d} and \eqref{rad} and will denote by $\mathcal{V}_k^\ast$ optimal coefficients defining a maximizer of $\sigma_k.$ The numerical maximization of the functional $\mathcal{C}_k$ may be difficult, in the sense that it is expected that the optimizers may have eigenvalues with large multiplicities,
\[\mathcal{C}_k(\mathcal{V}_k^\ast)=\mathcal{C}_{k+1}(\mathcal{V}_k^\ast)=...=\mathcal{C}_{k+\mathcal{M}-1}(\mathcal{V}_k^\ast),\]
for some $\mathcal{M}\geq1$ defining the multiplicity of the optimal eigenvalue. 

Thus, we must solve highly non-smooth optimizations. Several numerical approaches were proposed in the literature to deal with non-smoothness of the objective function in the optimization of eigenvalues (eg.~\cite{oudet,AntFr,AO,ak-kao-osting}). Here, we propose to use a gradient type method, with a convenient choice of the direction for performing a line search. We will assume that, due to numerical errors, all the domains that considered in the optimization procedure have simple eigenvalues and thus, we can calculate the gradients $g_k,g_{k+1},...,g_{k+\mathcal{M}-1}$ corresponding respectively to $\mathcal{C}_k,\mathcal{C}_{k+1},...,\mathcal{C}_{k+\mathcal{M}-1}$.

A typical situation during the optimization procedure is to obtain a certain domain having a few eigenvalues which are very close to each other and we would like the algorithm to increase all of them. To be more precise, let's assume that we define a small threshold parameter $\epsilon$ and that at some iteration, for some $\mathcal{M}\geq1,$ we obtain
\[\mathcal{C}_{k+\mathcal{M}}(\mathcal{V})-\mathcal{C}_k(\mathcal{V})\leq\epsilon,\ \text{but}\ \mathcal{C}_{k+\mathcal{M}+1}(\mathcal{V})-\mathcal{C}_k(\mathcal{V})>\epsilon\]
The ascent direction will be determined by solving the max min problem
\begin{equation}
\label{minmax}
\hat{v}=\max_{v\in\mathbb{R}^P:\left\|v\right\|=1}\min\left(g_{k}\cdot v,g_{k+1}\cdot v,...,g_{k+\mathcal{M}-1}\cdot v\right).\end{equation}
Note that from the computational point of view, the numerical solution of problem \eqref{minmax} is not expensive, when compared to the calculation of the eigenvalues and gradients.

\section{Numerical results and discussion}
\label{sec:numres}
In this section we present some numerical results that we obtained for the solution of Problem~\ref{shopteigvalprobnorm}. In all the experiments, we took N=20 for the parametrization of 3D and 4D domains, through functions $r$ defined in \eqref{linqcomb3d} and \eqref{rad} and the Steklov eigenvalues and eigenfunctions were calculated with 2000 and 8000 collocation points, respectively in 3D and 4D, taking $\delta=0.2$ in both cases.

\begin{figure}
 \centering
  \begin{tabular}{ccc}
 \includegraphics[width=0.3\textwidth]{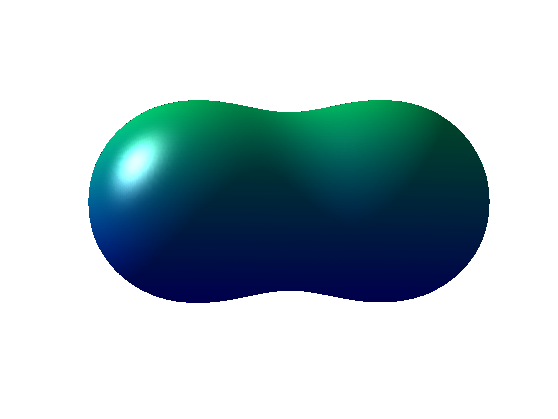}&
 \includegraphics[width=0.32\textwidth]{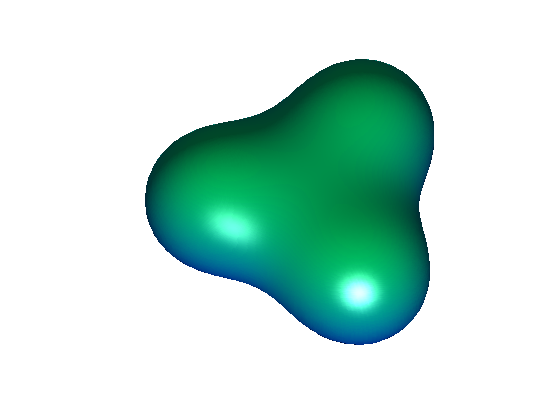}&
  \includegraphics[width=0.32\textwidth]{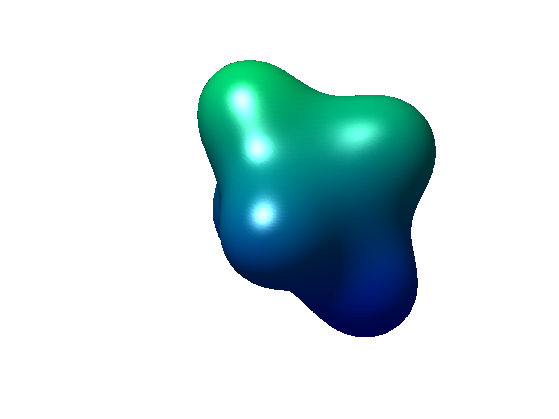}\\
  
    $\sigma^\star_2=2.19$ & $\sigma^\star_3=2.71$ & $\sigma^\star_5=3.58$ \\
  
   \includegraphics[width=0.32\textwidth]{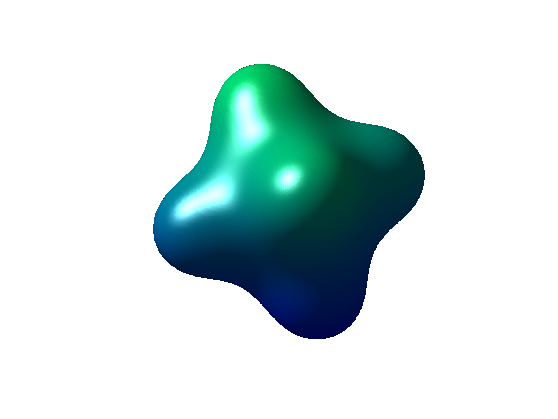}&
   \includegraphics[width=0.32\textwidth]{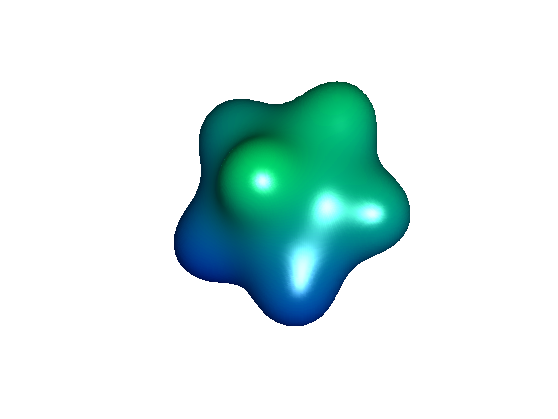}&
    \includegraphics[width=0.32\textwidth]{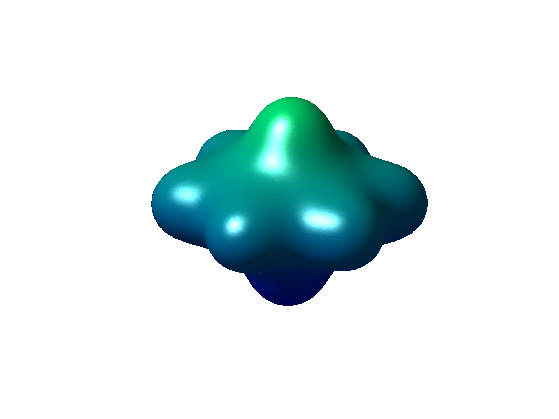}\\
    
        $\sigma^\star_6=3.97 $ & $\sigma^\star_7= 4.31$ & $\sigma^\star_8= 4.59$ \\
    
     \includegraphics[width=0.32\textwidth]{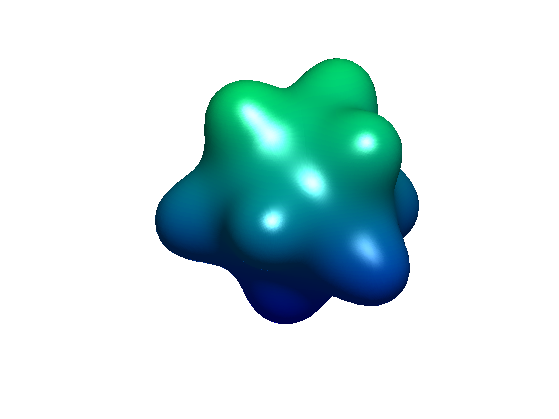}&
     \includegraphics[width=0.32\textwidth]{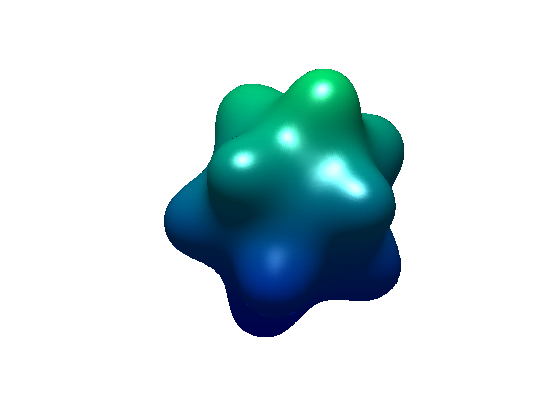}&
          \includegraphics[width=0.32\textwidth]{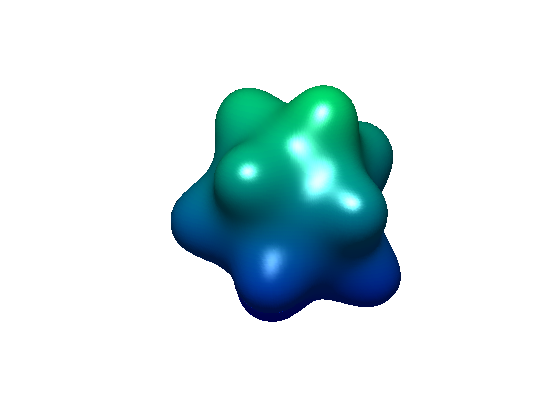}\\
     $\sigma^\star_{9}=4.90$ &      $\sigma^\star_{10}=5.17 $ & $\sigma^\star_{11}=5.44$  \\
             
   \includegraphics[width=0.32\textwidth]{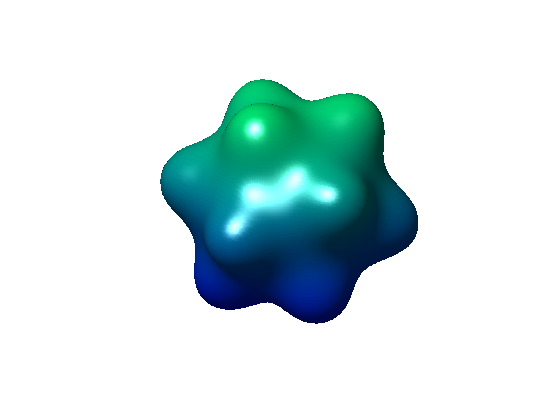}&
   \includegraphics[width=0.32\textwidth]{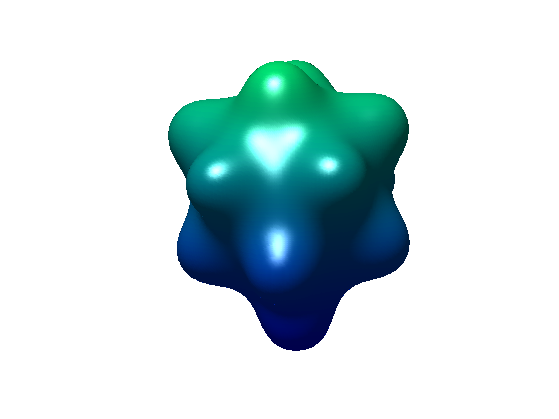}&
    \includegraphics[width=0.32\textwidth]{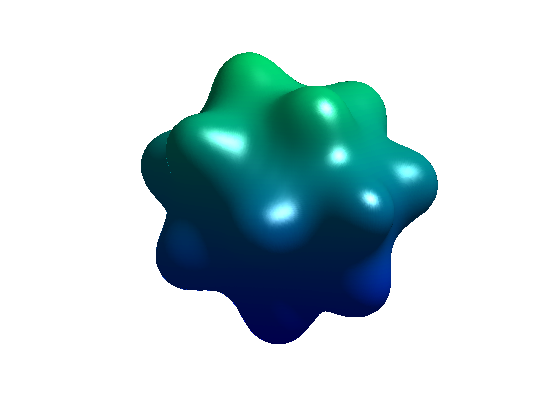}\\
    
  $\sigma^\star_{12}=5.76$ &   $\sigma^\star_{13}=5.96 $ & $\sigma^\star_{14}=6.20$   \\
        
    \includegraphics[width=0.32\textwidth]{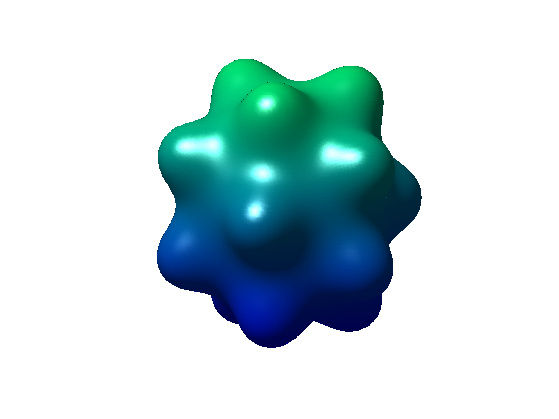}&
   \includegraphics[width=0.32\textwidth]{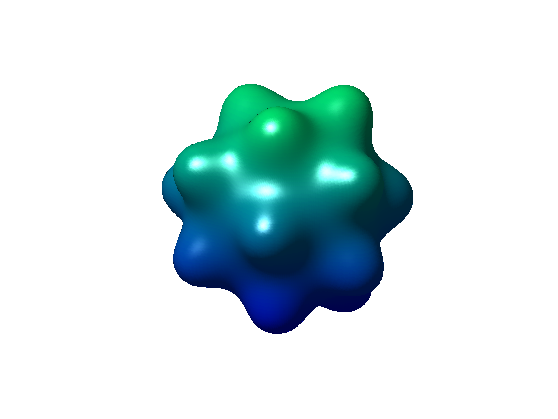}&
    \includegraphics[width=0.32\textwidth]{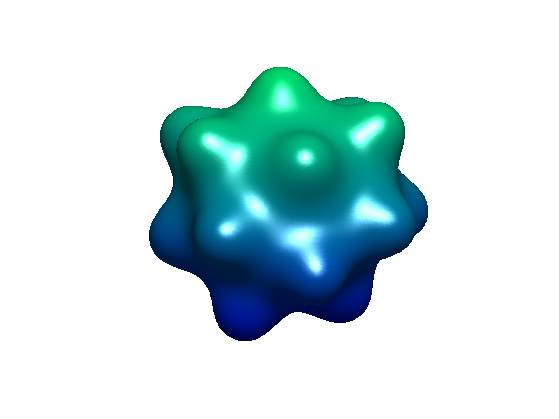}\\
    
  $\sigma^\star_{15}= 6.43$ &    $\sigma^\star_{16}=  6.65 $ & $\sigma^\star_{17}=6.84$  \\
        
     \includegraphics[width=0.32\textwidth]{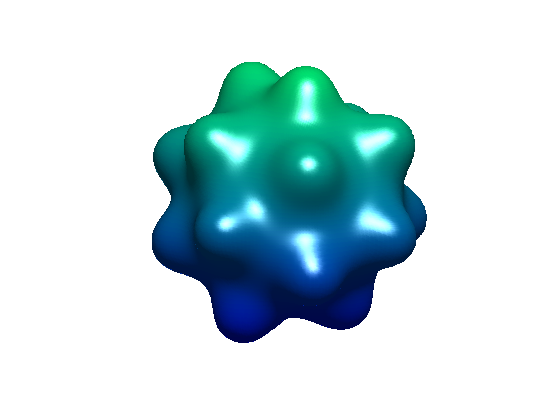}&
      \includegraphics[width=0.32\textwidth]{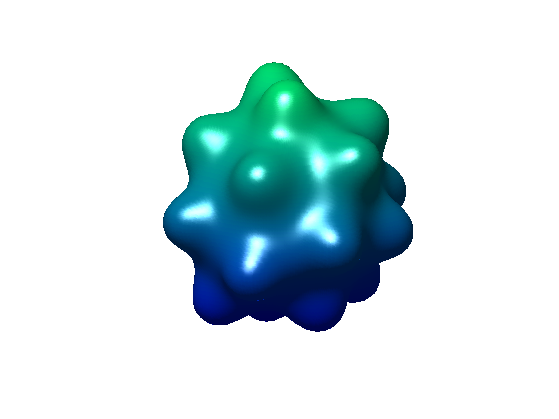}&
   \includegraphics[width=0.32\textwidth]{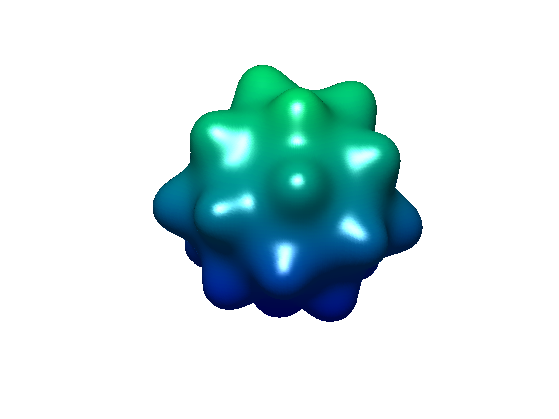}
   \\
    
 $\sigma^\star_{18}=7.06$  &   $\sigma^\star_{19}=7.25 $ & $\sigma^\star_{20}=7.46$  \\

 \end{tabular}
 \caption{Numerical optimizers and optimal eigenvalues obtained for Problem~\ref{shopteigvalprobnorm} in 3D. We obtained also the following optimal eigenvalue $\sigma^\star_4=3.22$, whose optimizer is the ball.}
 \label{fig:optimizers3d}
 \end{figure}

\begin{figure}
 \centering
  \begin{tabular}{ccc}
 \includegraphics[width=0.3\textwidth]{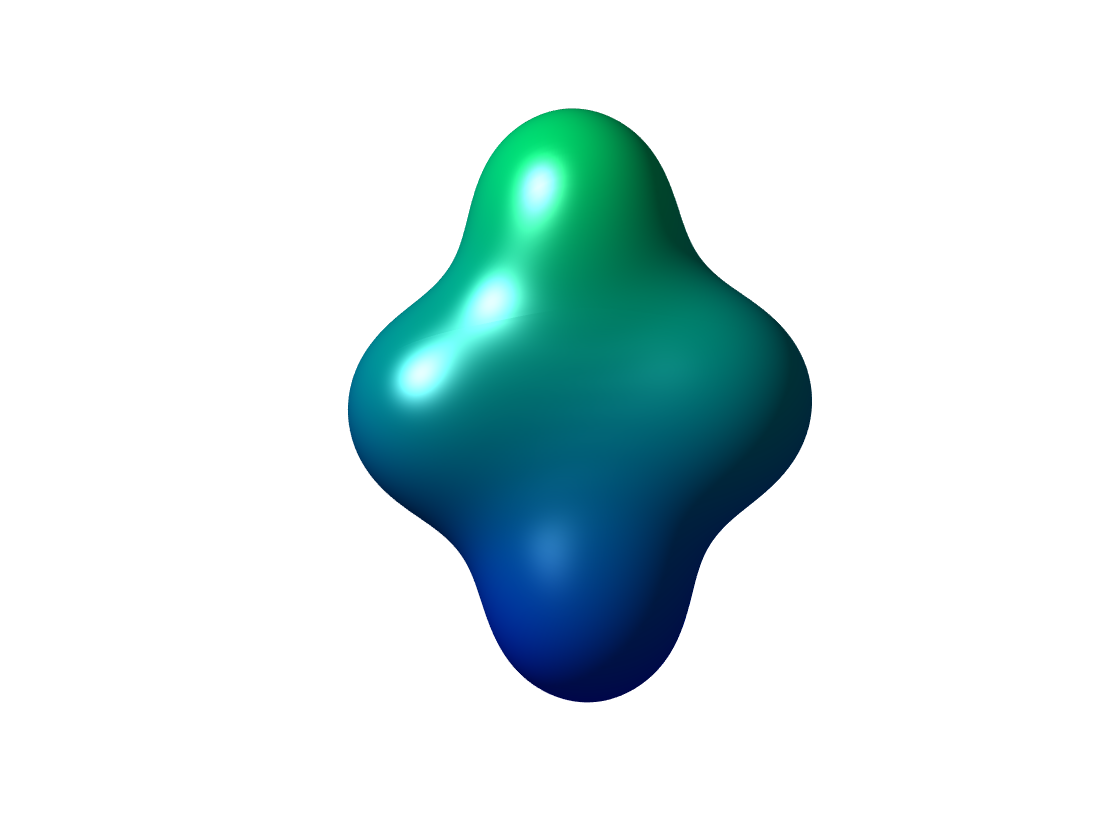}&
                               
 \end{tabular}
 \caption{Local maximizer of $\sigma_4$.}
 \label{fig:optimizers3dlocal}
 \end{figure}

\begin{figure}
 \centering
  \begin{tabular}{ccc}
 \includegraphics[width=0.4\textwidth]{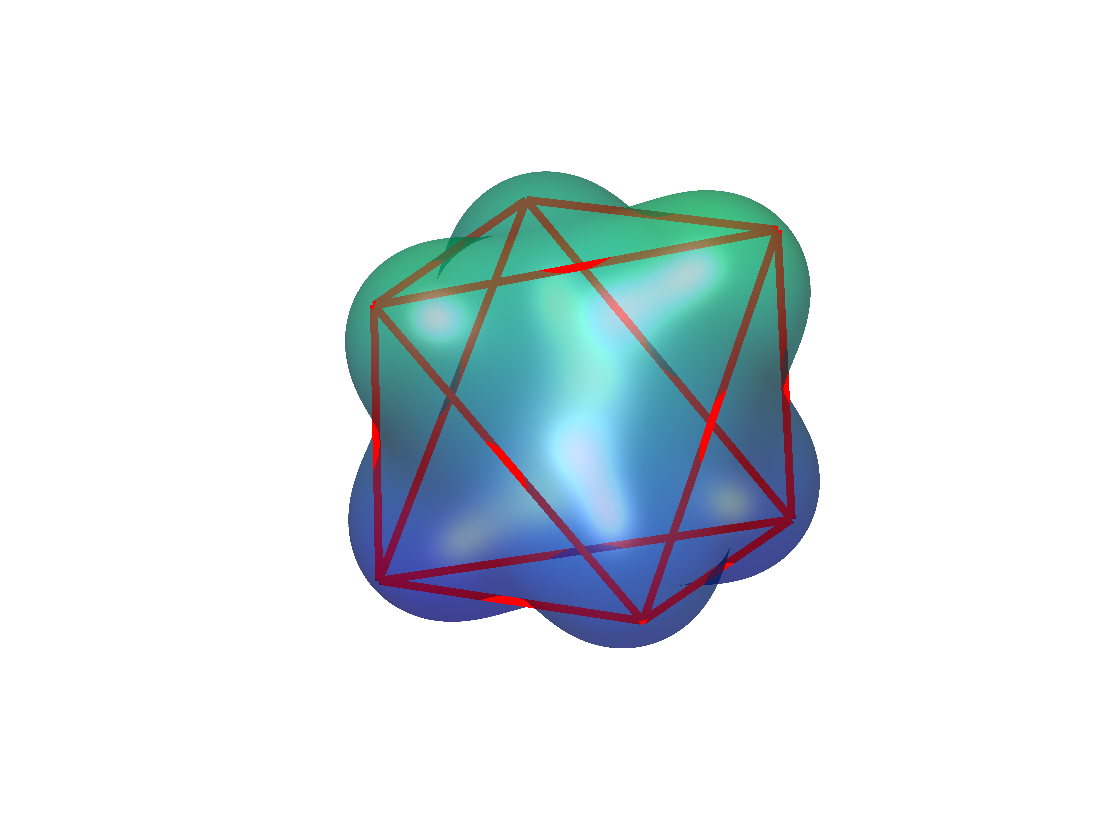}&
                               
 \end{tabular}
 \caption{Maximizer of $\sigma_6$ and an octahedron.}
 \label{fig:octaedro}
 \end{figure}

Figure~\ref{fig:optimizers3d} shows optimizers of Problem~\ref{shopteigvalprobnorm} among three dimensional geometries for $k=2,3,...,20$. The optimal eigenvalue is also marked in the Figure. All the values that are presented were obtained rounding down the optimal value obtained with our algorithm and are thus lower bounds for the optimal value. We can observe that, in a similar way that was obtained in the planar case in~\cite{ak-kao-osting}, in general the optimizers are well structured and have an increasing number of $\emph{buds}$. Indeed, the maximizer of $\sigma_k$ seems to have $k$ buds. However, there are some exceptions. We obtained a local maximizer for $\sigma_4,$ that is plotted in Figure~\ref{fig:optimizers3dlocal}, for which we obtain $\sigma_4\approx3.08,$ but this eigenvalue is smaller than the corresponding eigenvalue of the ball. Actually, our numerical results suggest that the optimizer of $\sigma_4$ is the ball. Some of optimizers seem to have some symmetries. For instance the optimizer of $\sigma_6$ seems to have the same symmetries of the octahedron, as illustrated in Figure~\ref{fig:octaedro}.

 \begin{table}
   \begin{tabular}{|c|c|}
   \hline
  Optimal eigenvalue & Multiplicity\\ 
  \hline\hline
  
 $\sigma^\star_2=2.19$ &3\\ 
\hline  $\sigma^\star_3=2.71$&4\\ 
\hline  $\sigma^\star_4=3.22$  &5\\ 
\hline  $\sigma^\star_5=3.58$&5\\ 
 \hline  $\sigma^\star_6=3.97 $  &6\\ 
\hline  $\sigma^\star_7= 4.31$&7\\ 
\hline   $\sigma^\star_8= 4.59$  &6\\ 
\hline      $\sigma^\star_{9}=4.90$  &6\\ 
\hline  $\sigma^\star_{10}=5.17 $  &7\\ 
\hline   $\sigma^\star_{11}=5.44$&7\\ 
 \hline   $\sigma^\star_{12}=5.76$ &9\\ 
\hline $\sigma^\star_{13}=5.96 $ &8\\ 
\hline  $\sigma^\star_{14}=6.20$  &7\\ 
 \hline   $\sigma^\star_{15}= 6.43$ &9\\ 
\hline $\sigma^\star_{16}=  6.65 $     &8\\ 
\hline $\sigma^\star_{17}=6.84$ &7\\ 
\hline   $\sigma^\star_{18}=7.06$   &8\\ 
\hline  $\sigma^\star_{19}=7.25 $ &7\\ 
 \hline $\sigma^\star_{20}=7.46$    &7\\
  \hline     
 \end{tabular}
 \caption{Optimal eigenvalues and multiplicities in 3D.}
 \label{tab:optimizers3d}
 \end{table}

\begin{figure}[ht]
\centering \includegraphics[width=0.32\textwidth]{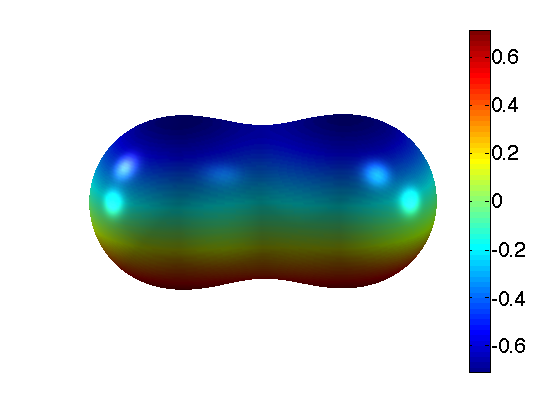}
\centering \includegraphics[width=0.32\textwidth]{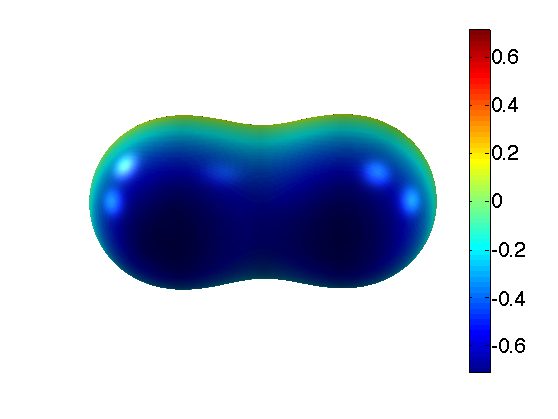}
\centering \includegraphics[width=0.32\textwidth]{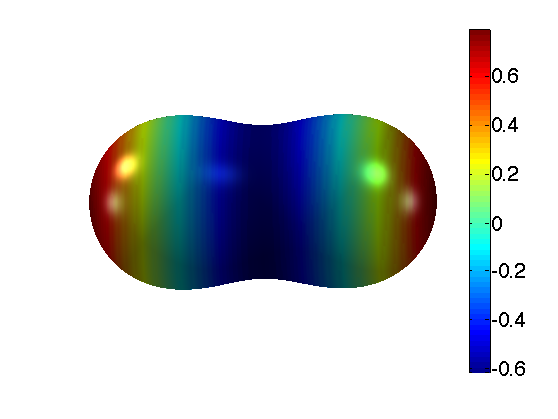}
\caption{Eigenfunctions associated to the optimal eigenvalue $\sigma_2$.} \label{fig:fplambda2}
\end{figure}

\begin{figure}[ht]
\centering \includegraphics[width=0.32\textwidth]{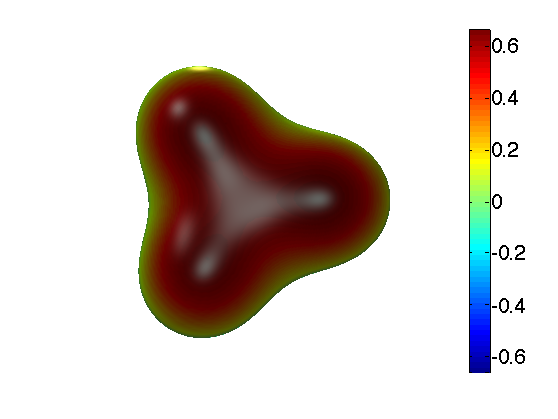}
\centering \includegraphics[width=0.32\textwidth]{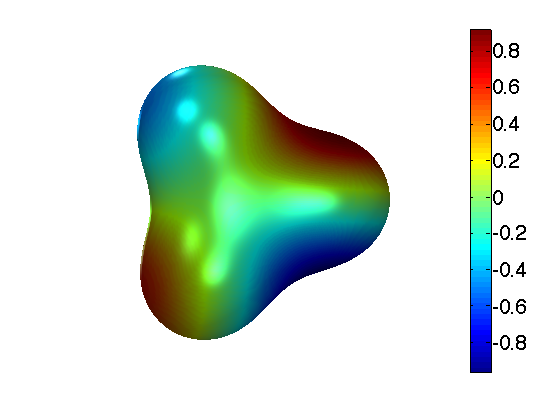}
\centering \includegraphics[width=0.32\textwidth]{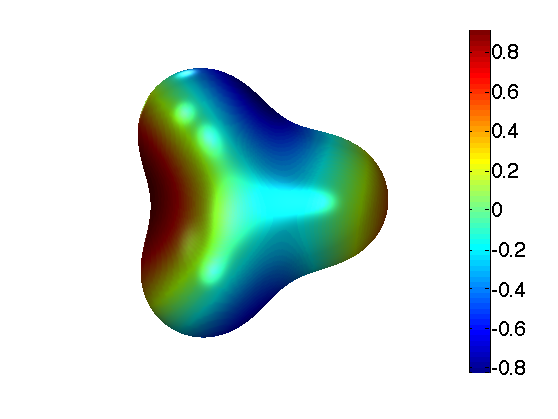}
\centering \includegraphics[width=0.32\textwidth]{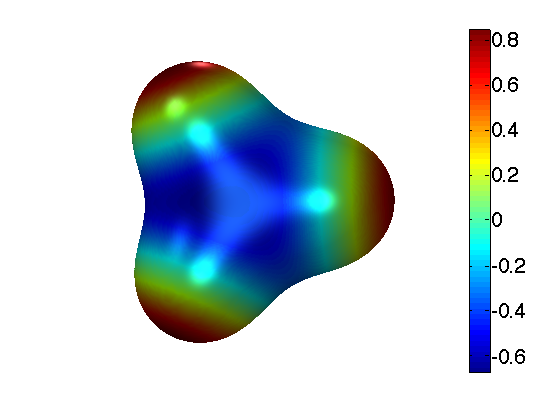}
\caption{Eigenfunctions associated to the optimal eigenvalue $\sigma_3$.} \label{fig:fplambda3}
\end{figure}

All the optimizers that we obtained have multiple optimal eigenvalues. For example the optimizer of $\sigma_2$ and $\sigma_3$ correspond to optimal eigenvalues with multiplicity three and four, respectively. In Figures~\ref{fig:fplambda2} and~\ref{fig:fplambda3} we plot  linear independent eigenfunctions associated to the optimal eigenvalue in both cases.

\begin{figure}[ht]
\centering \includegraphics[width=0.9\textwidth]{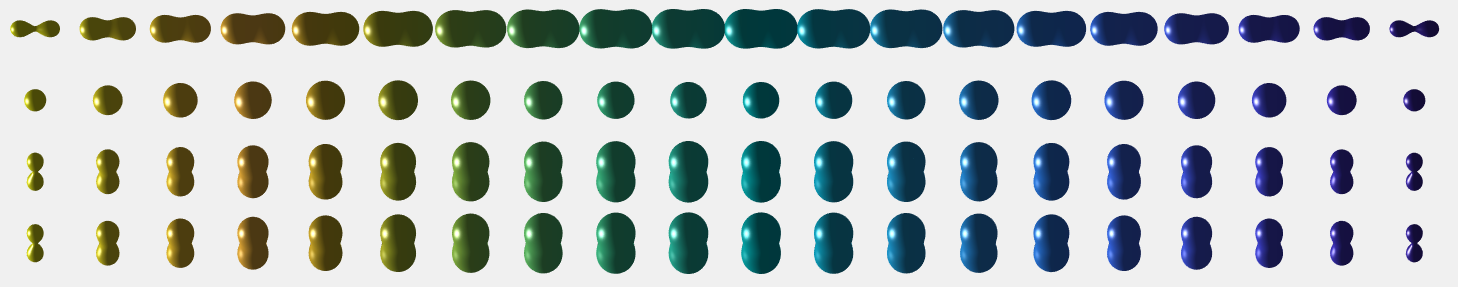}\\

\vspace{0.3cm}

\centering $\sigma_2^\ast=1.89$

\vspace{0.3cm}

\centering \includegraphics[width=0.9\textwidth]{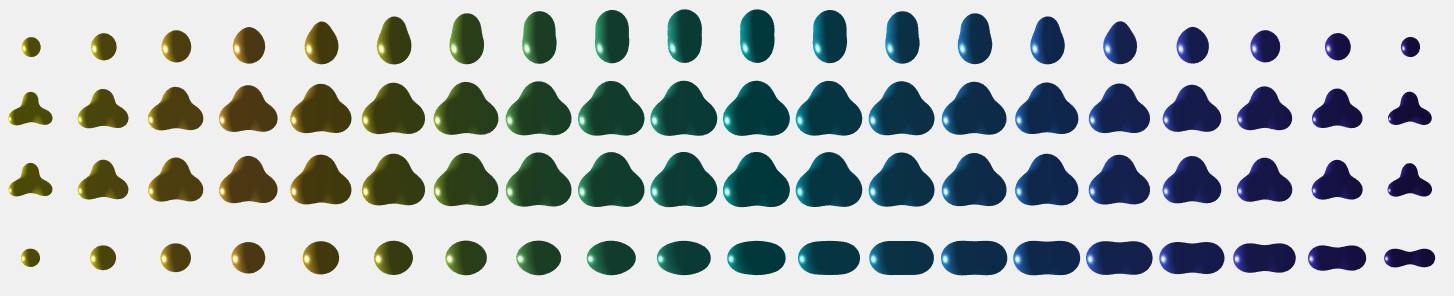}
\\

\vspace{0.3cm}

\centering $\sigma_3^\ast= 2.20$

\vspace{0.3cm}

\centering \includegraphics[width=0.9\textwidth]{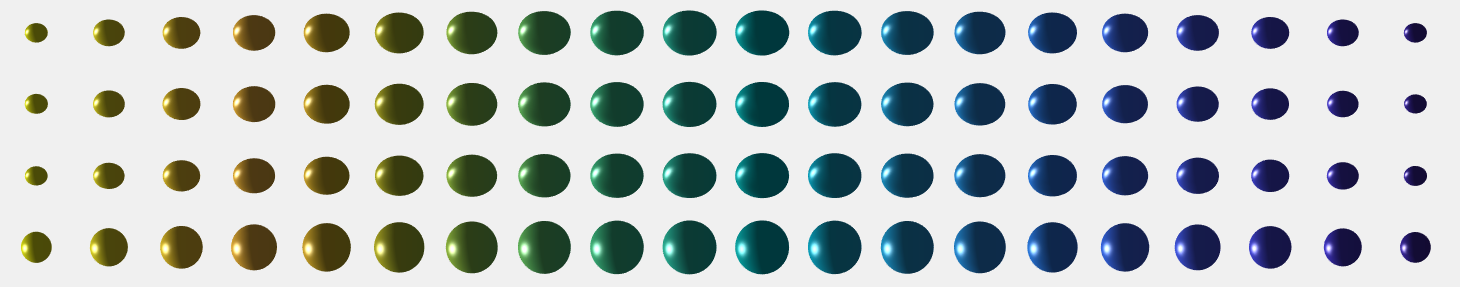}\\

\vspace{0.3cm}

\centering $\sigma_4^\ast=2.49$

\caption{Orthogonal cuts of the optimizers of $\sigma_k$, $k=2,3,4$ in 4D.} \label{fig:optimizers4d}
\end{figure}

Next, we show similar results obtained in 4D. Note that in this case is not trivial how to represent the optimizers. We used an algorithm proposed in~\cite{AO} that applies suitable rigid transformations to the optimizers obtained from the optimization procedure, in order to exhibit its symmetries (see~\cite{AO} for details). Figure~\ref{fig:optimizers4d} shows some 3D cuts of the optimizers in 4D in four orthogonal directions, for the optimizers of $\sigma_k$, $k=2,...,10$.  Every row in each picture corresponds to cuts in each of the four orthogonal directions. The optimal eigenvalues that were obtained are also presented at the Figure. The numerical results that we gathered suggest that the optimizer of $\sigma_5$ is the ball, for which we obtained $\sigma^\ast_5=3.01.$ Table~\ref{tab:optimizers4d} shows the optimal eigenvalues, together with the corresponding multiplicity. 

\begin{figure}[ht]

\centering \includegraphics[width=0.9\textwidth]{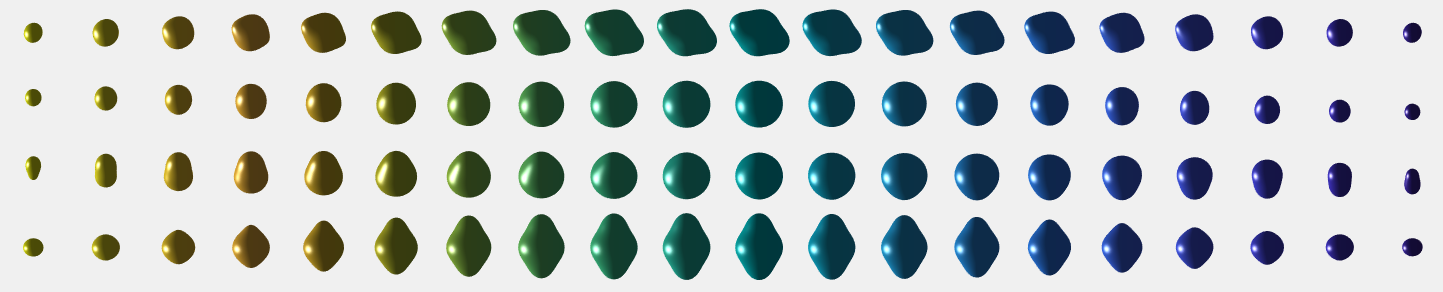}\\

\vspace{0.3cm}

\centering $\sigma_6^\ast= 3.12$

\vspace{0.3cm}

\centering \includegraphics[width=0.9\textwidth]{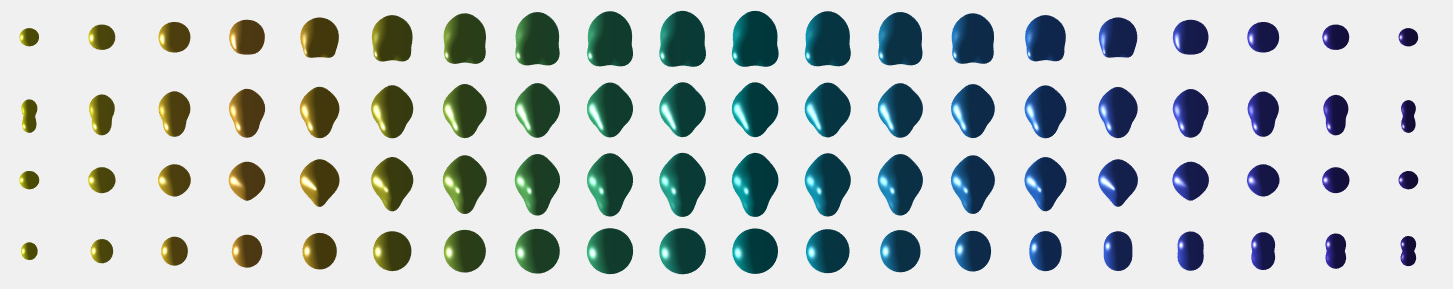}\\

\vspace{0.3cm}

\centering $\sigma_7^\ast= 3.23$

\vspace{0.3cm}

\centering \includegraphics[width=0.9\textwidth]{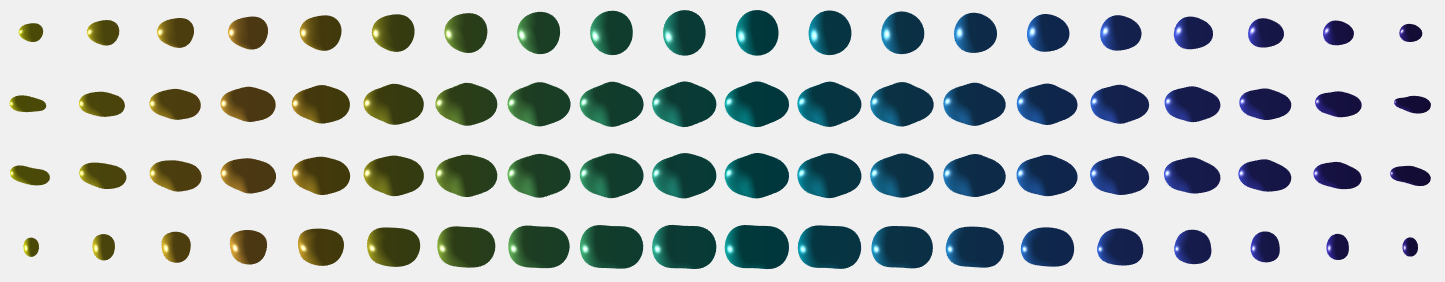}\\

\vspace{0.3cm}

\centering $\sigma_8^\ast= 3.36$

\vspace{0.3cm}

\centering \includegraphics[width=0.9\textwidth]{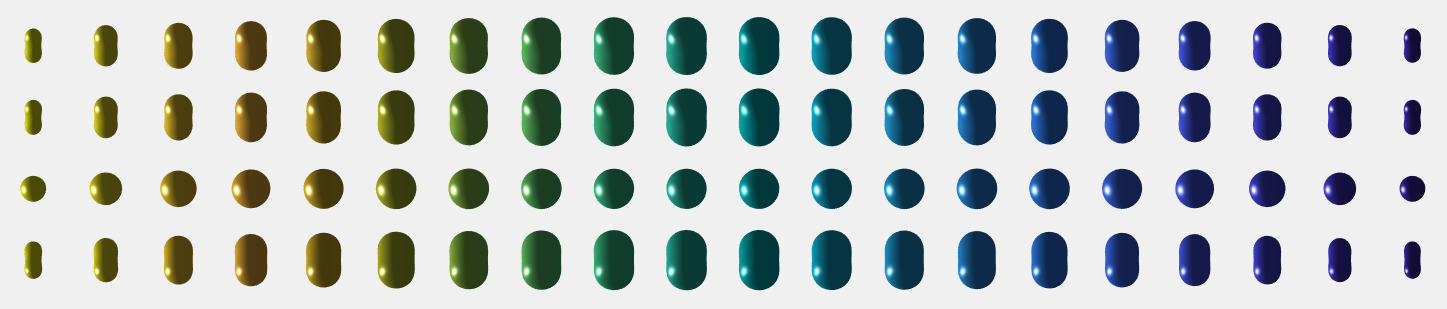}\\

\vspace{0.3cm}

\centering $\sigma_9^\ast= 3.60$

\vspace{0.3cm}

\centering \includegraphics[width=0.9\textwidth]{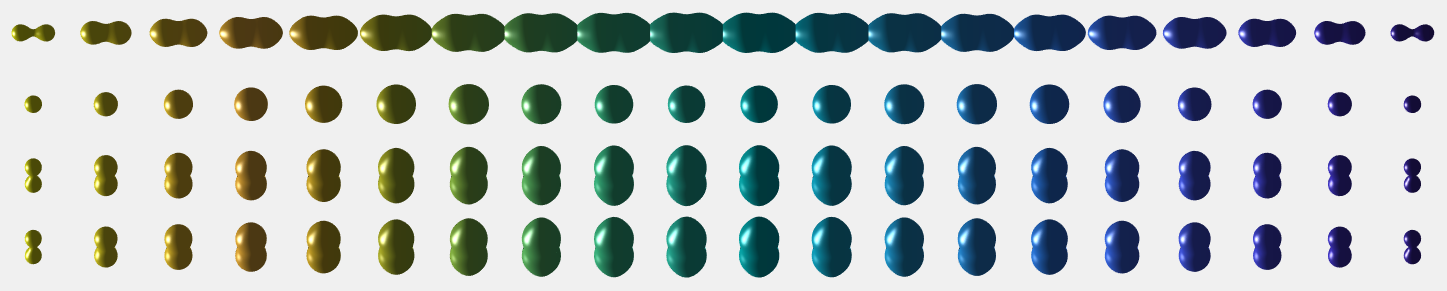}\\

\vspace{0.3cm}

\centering $\sigma_{10}^\ast= 3.71$

\vspace{0.3cm}

\caption{Orthogonal cuts of the optimizers of $\sigma_k$, $k=6,7,8,9,10$ in 4D.} \label{fig:optimizers4db}
\end{figure}

 \begin{table}
   \begin{tabular}{|c|c|}
   \hline
  Optimal eigenvalue & Multiplicity\\ 
  \hline\hline
  
 $\sigma^\star_2=1.89$ &4\\ 
\hline  $\sigma^\star_3=2.20$&5\\ 
\hline  $\sigma^\star_4=2.49$  &6\\ 
\hline  $\sigma^\star_5=3.01$&9\\ 
 \hline  $\sigma^\star_6=3.12 $  &9\\ 
\hline  $\sigma^\star_7= 3.23$&9\\ 
\hline   $\sigma^\star_8= 3.36$  &9\\ 
\hline      $\sigma^\star_{9}=3.60$  &9\\ 
\hline  $\sigma^\star_{10}=3.71 $  &9\\ 
  \hline     
 \end{tabular}
 \caption{Optimal eigenvalues and multiplicities in 4D.}
 \label{tab:optimizers4d}
 \end{table}

\end{document}